# About Three Dimensional Jump Boundary Value Problems for the Laplacian

Olexandr Polishchuk
Laboratory of Modeling and Optimization of Complex Systems
Pidstryhach Institute for Applied Problems of Mechanics and Mathematics, National Academy of Sciences of Ukraine
Lviv, Ukraine
od_polishchuk@ukr.net

*Abstract*—The conditions of well-posed solvability of searched function and its normal derivative three dimensional jump problem for the Laplacian and equivalent to them integral equation system for the sum of the simple and double layer potentials are determined in the Hilbert space, element of which as well as their normal derivatives have the jump through boundary surface.

## I. Introduction

Many physical processes (e.g. diffusion, heat flux, electrostatic field, perfect fluid flow, elastic motion of solid bodies, groundwater flow, etc.) are modeled using boundary value problems for Laplace equation [1, 2]. The powerful tools for solving such problems are potential theory methods, especially in the case of tired boundary surface or complex shape surface [3, 4]. These methods are a convenient for calculating desired solution in small domains [5]. In number of cases, application of potential theory methods requires determination the conditions of well-posed solvability for corresponding integral equations. Review of such conditions for main three dimensional boundary value problems for the Laplacian and equivalent to them integral equations for the simple and double layer potentials contains in [6]. These results allow us to use projection methods [7, 8] for numerical solution of such integral equations, avoiding the use of resource-consuming regularization procedures [9]. The need to determine the conditions of well-posed solvability also arises when the sum of simple and double layer potentials is used to solve the double-sided Dirichlet and Neumann problems [10] or double-sided Dirichlet-Neumann problem [11] in the space of functions that, same as their normal derivatives, have a jump on crossing boundary surface. When the domain and its environment have different physical properties, there is a need to solve the boundary value problems with jump conditions. Depending on the properties of searched solution when passing through the boundary surface, this can be a problem with the jump of searched function [12], a jump of its normal derivative [13], or a problem with both conditions simultaneously. The conditions of well-posed solvability of the latter problem in differential and integral formulations are investigated in this paper. Methods of integral equations and theory of boundary operators [6] enable not only to determine the properties of operators of such problems, but also to build effective methods for their numerical solution.

## II. Normal derivative jump problem

Let $G$ be the bounded open set in $R^3$, the boundary of which is Lipshitz surface $\Gamma$. Let us denote $G' = R^3 \setminus \overline{G}$ and introduce into $G, G'$ Sobolev spaces $W_2^1(G)$ and

$$W_{2,0}^1(G') = \{u^e \in D'(G') : u^e/r, Du^e \in L_2(G')\}.$$

Let us determine the Hilbert space $W_2^{1/2}(\Gamma)$ on the surface $\Gamma$ and introduce the space

$$H^1 = W_2^1(G) \times W_{2,0}^1(G').$$

Determine the linear continuous trace operators

$$\gamma_0^i u^i = u^i\big|_{\Gamma_i},\ \gamma_0^e u^e = u^e\big|_{\Gamma_e},\ \gamma_0^i : W_2^1(G) \to W_2^{1/2}(\Gamma),\ \gamma_0^e : W_{2,0}^1(G') \to W_2^{1/2}(\Gamma),$$

$$\gamma_1^i u^i = \partial u^i/\partial \mathbf{n}\big|_{\Gamma_i},\ \gamma_1^e u = \partial u/\partial \mathbf{n}\big|_{\Gamma_e},\ \gamma_1^i : W_2^1(G) \to W_2^{-1/2}(\Gamma),\ \gamma_1^e : W_{2,0}^1(G') \to W_2^{-1/2}(\Gamma)$$

where $\Gamma_i$ and $\Gamma_e$ are the internal and external sides of surface $\Gamma$, accordingly, $\mathbf{n}$ is a normal to the surface $\Gamma$ external to the domain $G$, and $W_2^{-1/2}(\Gamma)$ is the space duel to $W_2^{1/2}(\Gamma)$. Denote

$$\gamma_0^{\varepsilon_0} = \gamma_0^i - \varepsilon_0 \gamma_0^e,\ \gamma_1^{\varepsilon_1} = \gamma_1^i - \varepsilon_1 \gamma_1^e,\ \varepsilon_0, \varepsilon_1 > 0,$$

where $u = (u^i, u^e)$, $u^i \in W_2^1(G)$, $u^e \in W_{2,0}^1(G')$, and

$$W_2^1(G; \Delta = 0) = \{u^i \in W_2^1(G) : \Delta u^i = 0\},$$

$$W_{2,0}^1(G'; \Delta = 0) = \{u^e \in W_{2,0}^1(G') : \Delta u^e = 0\}.$$

For arbitrary $u^i \in W_2^1(G; \Delta = 0)$ and $v^i \in W_2^1(G)$ we have the Green's formula [12]

$$\int_G \nabla u^i \nabla v^i = \left\langle \gamma_1^i u^i, \gamma_0^i v^i \right\rangle_{W_2^{-1/2}(\Gamma) \times W_2^{1/2}(\Gamma)} \tag{1}$$

where $\langle .,. \rangle_{V \times V^*}$ is the duality relation on $V \times V^*$. For arbitrary $u^e \in W_{2,0}^1(G'; \Delta u = 0)$ and $v^e \in W_{2,0}^1(G')$ we have the Green's formula [14]

$$\int_{G'} \nabla u^e \nabla v^e = -\left\langle \gamma_1^e u^e, \gamma_0^e v^e \right\rangle_{W_2^{-1/2}(\Gamma) \times W_2^{1/2}(\Gamma)}. \tag{2}$$

Introduce the Hilbert space

$$H_\Gamma^1 = \{u \in H^1 : \gamma_0^1 u = 0\}, \quad (u,v)_{H_\Gamma^1} = \int_{G \cup G'} \nabla u \nabla v, \quad \|u\|_{H_\Gamma^1} = (u,u)_{H_\Gamma^1}^{1/2},$$

where $(.,.)_V$ is the scalar product on $V$. Introduce the space

$$H_{\Gamma, \Delta = 0}^1 = \{u \in H_\Gamma^1 : \Delta u = 0\}.$$

For arbitrary $u \in H_{\Gamma, \Delta = 0}^1$ and $v \in H_\Gamma^1$ from (1) and (2), we have the Green's formula

$$(u,v)_{H_\Gamma^1} = \left\langle \gamma_1^1 u, \gamma_0 v \right\rangle_{W_2^{-1/2}(\Gamma) \times W_2^{1/2}(\Gamma)} \tag{3}$$

where $\gamma_0 v = \gamma_0^i v^i = \gamma_0^e v^e$.

The next result is in order [12].

Theorem 1. Operator $\gamma_1^1$ is an isomorphism of $H_{\Gamma, \Delta = 0}^1$ onto $W_2^{-1/2}(\Gamma)$, and arbitrary function $u \in H_{\Gamma, \Delta = 0}^1$ can be represented as

$$u(x) = (U \gamma_1^1 u)(x) \equiv \frac{1}{4\pi} \int_\Gamma \frac{\gamma_1^1 u(y)}{|x - y|} d\Gamma_y, \quad x \in G, G', \quad y \in \Gamma.$$

Consider the boundary value problem: to find function

$$u \in H_{\Gamma, \Delta = 0}^1, \tag{4}$$

which satisfied condition

$$\gamma_1^{\varepsilon_1} u = f_1, \quad f_1 \in W_2^{-1/2}(\Gamma). \tag{5}$$

The next result is in order [12]

Theorem 2. Operator $\gamma_1^{\varepsilon_1}$ is an isomorphism of $H_{\Gamma, \Delta = 0}^1$ onto $W_2^{-1/2}(\Gamma)$ and equivalent to the problem (4)-(5) integral equation for the simple layer potential

$$(S\sigma)(x) \equiv (\gamma_1^{\varepsilon_1} U \sigma)(x) = f_1(x), \quad \sigma = \gamma_1^1 u \in W_2^{-1/2}(\Gamma), \quad f_1 \in W_2^{-1/2}(\Gamma), \tag{6}$$

has one and only one solution.

III. SEARCHED FUNCTION JUMP PROBLEM

Introduce the Hilbert space [14]

$$K_\Gamma^1 = \{u \in H^1 \setminus R : \gamma_1^1 u = 0\}$$

where $R$ is the space of functions that are constant on $G$ and vanishing on $G'$, and

$$(u,v)_{K_\Gamma^1} = \int_{G\cup G'} \nabla u \nabla v, \quad \|u\|_{K_\Gamma^1} = (u,u)_{K_\Gamma^1}^{1/2}.$$

Introduce the spaces

$$K_{\Gamma,\Delta=0}^1 = \{u \in K_\Gamma^1 : \Delta u = 0\},$$

$$\widetilde{W}_2^{1/2}(\Gamma) = W_2^{1/2}(\Gamma) \setminus P,$$

$$\widetilde{W}_2^{-1/2}(\Gamma) = \{g \in W_2^{-1/2}(\Gamma) : \langle g, 1 \rangle_{W_2^{1/2}(\Gamma) \times W_2^{-1/2}(\Gamma)} = 0\}$$

where $P$ is the set of constant functions on $\Gamma$. For arbitrary $u \in K_{\Gamma,\Delta=0}^1$ and $v \in K_\Gamma^1$ from (1) and (2), we have the Green's formula

$$(u,v)_{K_\Gamma^1} = \langle \gamma_0^1 u, \gamma_1 v \rangle_{W_2^{1/2}(\Gamma) \times W_2^{-1/2}(\Gamma)} \tag{7}$$

where $\gamma_1 v = \gamma_1^i v^i = \gamma_1^e v^e$.

The next result is in order [14].

Theorem 3. Operator $\gamma_0^1$ is an isomorphism of $K_{\Gamma,\Delta=0}^1$ onto $\widetilde{W}_2^{1/2}(\Gamma)$, and arbitrary function $v \in K_{\Gamma,\Delta=0}^1$ can be represented as

$$v(x) = (V\gamma_0^1 v)(x) \equiv -\frac{1}{4\pi} \int_\Gamma \gamma_0^1 v(y) \frac{\partial}{\partial \mathbf{n}_y} \frac{1}{|x-y|} d\Gamma_y, \quad x \in G, G', \quad y \in \Gamma.$$

Consider the boundary value problem: to find function

$$v \in K_{\Gamma,\Delta=0}^1, \tag{8}$$

which satisfied condition

$$\gamma_0^{\varepsilon_0} v = f_0, \quad f_0 \in \widetilde{W}_2^{1/2}(\Gamma). \tag{9}$$

The next result is in order [13]

Theorem 4. Operator $\gamma_0^{\varepsilon_0}$ is an isomorphism of $K_{\Gamma,\Delta=0}^1$ onto $\widetilde{W}_2^{1/2}(\Gamma)$ and equivalent to the problem (8)-(9) integral equation for the double layer potential

$$(Dq)(x) \equiv (\gamma_0^{\varepsilon_0} V q)(x) = f_0(x), \tag{10}$$

$$q = \gamma_1^0 u \in \widetilde{W}_2^{1/2}(\Gamma), \quad f_0 \in \widetilde{W}_2^{1/2}(\Gamma),$$

has one and only one solution.

IV. SEARCHED FUNCTION AND ITS NORMAL DERIVATIVE JUMP PROBLEM

Spaces $H_{\Gamma,\Delta=0}^1$ and $K_{\Gamma,\Delta=0}^1$ are ortogonal relatively introduced scalar product [15]. Introduce the space

$$HK_\Gamma^1 = H_\Gamma^1 \oplus K_\Gamma^1, \quad (u,v)_{HK_\Gamma^1} = \int_{G\cup G'} \nabla u \nabla v, \quad \|u\|_{HK_\Gamma^1} = (u,u)_{HK_\Gamma^1}^{1/2}$$

and

$$HK_{\Gamma,\Delta=0}^1 = H_{\Gamma,\Delta=0}^1 \oplus K_{\Gamma,\Delta=0}^1.$$

From (1) and (2) we obtain the Green's formula

$$(u,v)_{H_\Gamma^1 \cup K_\Gamma^1} = \langle \gamma_1^i u, \gamma_0^i v \rangle_{W_2^{-1/2}(\Gamma) \times W_2^{1/2}(\Gamma)} - \langle \gamma_1^e u, \gamma_0^e v \rangle_{W_2^{-1/2}(\Gamma) \times W_2^{1/2}(\Gamma)} \tag{11}$$

for arbitrary $u \in HK_{\Gamma,\Delta=0}^1$ and $v \in HK_\Gamma^1$.

Then from theorems 2 and 4 we have that the following result is in order [15].

Theorem 5. Operator $(\gamma_1^1, \gamma_0^1)$ is an isomorphism of $HK_{\Gamma,\Delta=0}^1$ onto $W_2^{-1/2}(\Gamma) \times \widetilde{W}_2^{1/2}(\Gamma)$ and arbitrary function $w \in HK_{\Gamma,\Delta=0}^1$ can be singlevaluedly represent as

$$w(x) = (W(\gamma_1^1 u, \gamma_0^1 v))(x) = (U, V) \, (\gamma_1^1 u, \gamma_0^1 v)^T (x) = = (U \, \gamma_1^1 u)(x) + (V \, \gamma_0^1 v)(x) \,, \tag{12}$$

where $u \in H^1_{\Gamma, \Delta=0}$, $v \in K^1_{\Gamma, \Delta=0}$, $x \in G, G'$.

Consider the boundary value problem: to find function

$$w \in HK^1_{\Gamma, \Delta=0}, \tag{13}$$

which satisfied conditions

$$\gamma_0^{\varepsilon_0} w = g_0, \; g_0 \in \widetilde{W}_2^{1/2}(\Gamma), \tag{14}$$

$$\gamma_1^{\varepsilon_1} w = g_1, \; g_1 \in W_2^{-1/2}(\Gamma). \tag{15}$$

Taking into account representation (12) and denote

$$\gamma_0^i u^i = \gamma_0^e u^e = \gamma_0 u$$

for $u \in H^1_{\Gamma, \Delta=0}$ and

$$\gamma_1^i v^i = \gamma_1^e v^e = \gamma_1 v$$

for $v \in K^1_{\Gamma, \Delta=0}$, from (14)-(15) obtain

$$(1 - \varepsilon_0) \gamma_0 u + \gamma_0^{\varepsilon_0} v = g_0, \tag{16}$$

$$\gamma_1^{\varepsilon_1} u + (1 - \varepsilon_1) \gamma_1^{\varepsilon_1} v = g_1. \tag{17}$$

Determine

$$\gamma_0 u = \Phi \gamma_1^{\varepsilon_1} u,$$

where $\Phi$ is a boundary operator which matches to the value of jump of the normal derivative of searched function its trace on the surface $\Gamma$. Operator $\Phi$, constructed in the form

$$\Phi = \widetilde{S} \, S^{-1},$$

where $\widetilde{S} = \gamma_0 U$ and $S$ is determined according to (6) is an isomorphism of $W_2^{-1/2}(\Gamma)$ onto $W_2^{1/2}(\Gamma)$ [6]. Note that $\widetilde{S}$ and $S$ are operators of integral equations for the simple layer potential equivalent to the Dirichlet and searched function normal derivative jump problems respectively. Algorithmically, operator $\Phi$ is realized by solving the integral equation for simple layer potential equivalent to the normal derivative jump problem and calculating by the finite formulas the values of searched function on the boundary surface as a trace of simple layer potential with the density determined in the previous step.

Determine

$$\gamma_1 v = \Psi \gamma_0^{\varepsilon_0} v,$$

where $\Psi$ is a boundary operator which matches to the value of searched function jump its trace on the surface $\Gamma$. Operator $\Psi$, constructed in the form

$$\Psi = \widetilde{D} \, D^{-1},$$

where $\widetilde{D} = \gamma_1 V$ and $D$ is determined according to (10) is an isomorphism of $\widetilde{W}_2^{1/2}(\Gamma)$ onto $\widetilde{W}_2^{-1/2}(\Gamma)$ [6]. Note that $\widetilde{D}$ and $D$ are operators of integral equations for the double layer potential equivalent to the Neumann and searched function jump problems respectively. Algorithmically, operator $\Psi$ is realized by solving the integral equation for double layer potential equivalent to the searched function jump problem and calculating by the finite formulas the values of its normal derivative on the boundary surface as a trace of double layer potential with the density determined in the previous step.

Then from (16) - (17) consistently obtain

$$(1 - \varepsilon_0) \, \Phi \gamma_1^{\varepsilon_1} u + \gamma_0^{\varepsilon_0} v = g_0,$$

$$\Psi^{-1} \gamma_1^{\varepsilon_1} u / (1 - \varepsilon_1) + \gamma_0^{\varepsilon_0} v = \Psi^{-1} g_1 / (1 - \varepsilon_1),$$

$$(1 - \varepsilon_0) \, \Phi \gamma_1^{\varepsilon_1} u - \Psi^{-1} \gamma_1^{\varepsilon_1} u / (1 - \varepsilon_1) = g_0 - \Psi^{-1} g_1 / (1 - \varepsilon_1),$$

$$[(1-\varepsilon_0)(1-\varepsilon_1)\Phi - \Psi^{-1}]\gamma_1^{\varepsilon_1} u = (1-\varepsilon_1)g_0 - \Psi^{-1}g_1,$$

$$\gamma_1^{\varepsilon_1} u = p_1, \qquad (18)$$

where

$$p_1 = [(1-\varepsilon_0)(1-\varepsilon_1)\Phi^{-1} - \Psi] \times [(1-\varepsilon_1)g_0 - \Psi^{-1}g_1)] \in W_2^{-1/2}(\Gamma).$$

From theorem 2 follows that problem (4), (18) has a unique solution. Solving this problem is equivalent to solving an integral equation for simple layer potential

$$(S\sigma)(x) \equiv (\gamma_1^{\varepsilon_1} U\sigma)(x) = p_1(x), \; \sigma \in W_2^{-1/2}(\Gamma), \; p_1 \in W_2^{-1/2}(\Gamma), \qquad (19)$$

which, by theorem 2, also has a unique solution.

Next, from (16) obtain

$$\gamma_0^{\varepsilon_0} v = g_0 - (1-\varepsilon_0)\gamma_0 u$$

and

$$\gamma_0^{\varepsilon_0} v = p_0, \qquad (20)$$

where

$$p_0 = g_0 - (1-\varepsilon_0)\Phi p_1.$$

From theorem 4 follows that if the condition

$$p_0 \in \widetilde{W}_2^{1/2}(\Gamma) \qquad (21)$$

holds, then problem (8), (20) has a unique solution. Solving this problem is equivalent to solving an integral equation for double layer potential

$$(Dq)(x) \equiv (\gamma_0^{\varepsilon_0} Vq)(x) = p_0(x), \; q \in \widetilde{W}_2^{1/2}(\Gamma), \; p_0 \in \widetilde{W}_2^{1/2}(\Gamma), \qquad (22)$$

which, by theorem 4, also has a unique solution.

Given the results of theorems 2 and 4, we obtain the validity of following statement.

Theorem 6. Under condition (21), the problem (13)-(15) reduces to the sequential solution of problems (4), (18) and (8), (20) or equivalent to them integral equations for the simple and double layers potentials (19) and (22), which have the unique solutions. Then the solution of problem (13)-(15) is determined by relation (12).

## V. Conclusions

Using for the approximation of unknown densities of the potentials a system of $N$ linearly independent functions (Lagrangian finite elements, B-splines, etc. [7, 8]), we arrive for the solution of problem (13) - (15) by means of the sum of the simple and double layer potentials the need to solve the system of linear algebraic equations (SLAE) with dense matrices of dimension $2N$, which requires the performing of $8O(N^3)$ operations. Using the procedure implemented to prove theorem 6 requires the solution of five SLAEs with matrices of dimension $N$, which requires the performing of $5O(N^3)$ operations, which is almost 40% less. Thus, applying the theory of boundary operators allows us not only to determine the conditions of well-posed solvability of separate boundary value problems for the Laplace equation, but also to build effective methods for their numerical solution using the simple and double layer potentials.